# THE SET OF INFINITE VALENCE VALUES OF AN ANALYTIC FUNCTION

Julian Gevirtz
2005 North Winthrop Road
Muncie, Indiana
USA
jgevirtz@gmail.com

## ABSTRACT

It is shown (Theorem A and its corollary) that if $g$ is any nonconstant nonunivalent analytic function on a half-plane $\mathbb{H}$ and if $D$ is either a half-plane or a smoothly bounded Jordan domain, then there is a function $f$ on $D$ for which $f'(D) \subset g'(\mathbb{H})$ such that for any neighborhood $U$ of any point of $f(\partial D)$ the set of values $w \in U$ which $f$ assumes infinitely many times in $D$ has Hausdorff dimension 1. From this it follows (Theorem C) that in the Becker univalence criteria for the disc and upper-half plane ($|f''(z)/f'(z)| \leq \frac{1}{1-|z|^2}$ and $|f''(z)/f'(z)| \leq \frac{1}{2\Im\{z\}}$, respectively) if the 1 in the numerator is replaced by any larger number, then there are functions $f$ satisfying the resulting bounds the set of whose infinitely assumed values has this same dimension 1 property.

## 1. Introduction.

In what follows we study the set of values of infinite valence of an analytic function $f$ on a domain $D$, that is, the set $\mathcal{I}(f)$ of points $w \in f(D)$ for which $f^{-1}(w)$ is infinite. We will limit the discussion to $D \in \mathfrak{D}$, the class consisting of the upper half-plane $\mathbb{H}$ together with all smoothly bounded Jordan domains. The word "smooth" is used here somewhat vaguely, $C^\infty$ being both clearly enough and clearly far too much. It is immediate that if $D$ is a bounded domain and $f$ has a continuous extension to $\overline{D}$, then $\mathcal{I}(f) \subset f(\partial D)$ since in this case, if $f^{-1}(w)$ is infinite then $f^{-1}(w)$ has an accumulation point $w_0 \in \partial D$. Let

$$R(m, M) = \{z : m < |z| < M\},$$
$$\Delta(a, r) = \{z : |z - a| < r\}, \Delta(r) = \Delta(0, r), \Delta = \Delta(1),$$

For $X \subset \mathbb{C}$, $\dim(X)$ and $\lambda_1(X)$ will denote, respectively, the Hausdorff dimension of $X$ and its 1-dimensional measure. We denote the symmetric difference $(X \backslash Y) \cup (Y \backslash X)$ by $X \Delta Y$. The main result to be established is the following

**Theorem A**. *Let $g$ be analytic and nonunivalent on $\mathbb{H}$, let $0 \notin g'_0(\mathbb{H})$ and let $D \in \mathfrak{D}$. Then there is an analytic function $\rho$ on $D$, for which $\rho(D) \subset \mathbb{H}$ such that for the antiderivatives $f$ of $g'_0(\rho(z))$, $\dim(\Delta(z_0, r) \cap \mathcal{I}(f)) = 1$ for all $r > 0$ and all $z_0 \in f(\partial D)$.*

The proof of this theorem, which we give in §3, is to some extent a refinement of that of [G4, Lemma 6, p.190] which says that *for every nonunivalent function $g$ on $\mathbb{H}$ and every $D$ which is either a smoothly bounded Jordan domain or a half-plane there is an infinitely valent $f$ on $D$ which satisfies $f'(D) \subset g'(\mathbb{H})$*. We mention the following corollary of Theorem A which is relevant to the general theory of first-order univalence criteria to be discussed in §4.

**Corollary**. *Let $g$ be analytic, nonconstant and nonunivalent on $\mathbb{H}$ and let $D \in \mathfrak{D}$. Then there is an analytic function $f$ on $D$ such that $f'(D) \subset g'(\mathbb{H})$ and $\dim(\Delta(z_0, r) \cap \mathcal{I}(f)) = 1$ for all $r > 0$ and all $z_0 \in f(\partial D)$.*

This follows immediately from Theorem A, since if $g'$ vanishes on $\mathbb{H}$, then $g'(\mathbb{H})$ contains $R(m, M)$ with arbitrarily large $M/m$ and since for sufficiently large $M/m$ there are nonunivalent functions on $\mathbb{H}$ the range of whose derivative lies in $R(m, M)$; indeed, such functions exist as long as $M/m$ is greater than the John constant of $\mathbb{H}$, which is known to be less than 7.19 (see [G1]). Theorem A is perhaps a little surprising in light of the following additional fact, whose simple proof we give in §2.

**Theorem B**. *Let $D$ be a smoothly bounded Jordan domain and let $f'(D) \subset R(m, M)$ for some $0 < m < M < \infty$. Then $\lambda_1(\mathcal{I}(f)) = 0$.*

Thus, at least in the case in which $g'(\mathbb{H}) \subset R(m, M)$, for any $f$ on $D$ for which $f'(D) \subset g'(\mathbb{H})$, $\mathcal{I}(f)$ must have 1-dimensional measure 0 but, on the other hand, Theorem A says that if $g$ is not univalent on $\mathbb{H}$, $\mathcal{I}(f)$ can actually be a dense 1-dimensional subset of $f(\partial D)$. In §4, in addition to discussing a sharp first-order univalence criterion which is a weak form of the Becker criterion,

$$|f''(z)/f'(z)| \leq 1/(1 - |z|^2) \Rightarrow f \text{ is univalent in } \Delta,$$

we apply Theorem A to show that

**Theorem C.** *Let $D = \Delta$ or $\mathbb{H}$. For any $\tau > 1$ there are functions on $D$ which satisfy*

$$|f''(z)/f'(z)| \leq \tfrac{\tau}{1-|z|^2}, \quad |f''(z)/f'(z)| \leq \tfrac{\tau}{2\Im\{z\}},$$

*for the cases $D = \Delta$, $D = \mathbb{H}$, respectively, such that the values of $f'$ all lie in some annulus $R(m, M)$ and such that for any $r > 0$ and $z_0 \in f(\partial D)$, $\dim(\Delta(z_0, r) \cap \mathcal{I}(f)) = 1$.*

This theorem is a considerable strengthening of results of Aksent'ev [A] and Avhadiev and Kayumov [AK] to the effect that for any $\tau > 1$ and any integer $k$ there are $k$-univalent functions $f_\tau$ on $\mathbb{H}$ and $\Delta$, which satisfiy $|f''_\tau(z)/f'_\tau(z)| \leq \tfrac{\tau}{2\Im\{z\}}$ and $|f''_\tau(z)/f'_\tau(z)| \leq \tfrac{\tau}{1-|z|^2}$, respectively. (Application of the above-mentioned lemma of [G4] already shows that there are $f_\tau$ satisfying these bounds which are infinitely valent.) We mention that for any nonunivalent $g$ on $\mathbb{H}$ a construction even simpler than that used to establish Theorem A yields an $f$ on $\mathbb{H}$ with $f'(\mathbb{H}) \subset g'(\mathbb{H})$ and for which $\lambda_1(\mathcal{I}(f) > 0$. This can be applied to the Becker criterion to show that if $\tau > 1$, then there is a function $f$ on $\mathbb{H}$ satisfying $|f''(z)/f'(z)| \leq \tfrac{\tau}{2\Im\{z\}}$ for which $\lambda_1(\mathcal{I}(f)) > 0$. (This difference between the case of $\mathbb{H}$ and that of Jordan domains arises because in the former context it is not necessarily true that if $f^{-1}(w)$ is infinite then $f^{-1}(w)$ has an accumulation point in $\partial \mathbb{H}$.)

## 2. Proof of Theorem B.

Let $A \subset \partial D$ be the set points $\zeta \in \partial D$ at which $f'(z)$ has a nontangential limit $f'(\zeta)$ as $z \to \zeta$ in $D$. As is well known, bounded analytic functions in $\Delta$ have nontangential limits at almost all points of $\partial \Delta$, so that the corresponding fact also holds for all $D \in \mathfrak{D}$. In examining what happens in the vicinity of a point $\zeta \in A$ we may, for the sake of simplicity and without loss of generality, assume that $\zeta = 0$ and that the positive tangent direction on $\partial D$ at $\zeta$ is that of the positive real axis. Let $\epsilon > 0$ and let $W(\epsilon)$ denote the sector $\{z \in \mathbb{H} : \epsilon < \arg\{z\} < \pi - \epsilon\}$. In light the existence of the nontangential limit of $f'$ at $\zeta$, there is a $\delta_0 = \delta_0(\epsilon) > 0$ such that for all positive $\delta < \delta_0$

(i) $\Delta(\delta) \cap W(\epsilon) \subset D$

(ii) $|f'(z) - f'(0)| < \epsilon$ for $z \in \Delta(\delta) \cap W(\epsilon)$.

In addition, simple trigonometry and the smoothness of $\partial D$ imply that for sufficiently small $\delta$

(iii) for every point $p \in D \cap \Delta(\delta)$ there is a $q \in \Delta(\delta) \cap W(\epsilon)$ for which the open segment $(p, q) \subset D$ and $|q - p| < 2\epsilon |p|$.

Let $p \in D \cap \Delta(\delta)$ and let $q$ be as in (iii). By (i) and the fact that $f'(0) \in \overline{R(m, M)}$ we have from (ii) that $|f(q) - f(0)| \geq (m - \epsilon)|q|$. On the other hand, from (iii) it follows that $|f(q) - f(p)| < 2\epsilon M|p|$, so that for all $p \in D \cap \Delta(\delta)$ we have

$$\begin{aligned} |f(p) - f(0)| &= |f(p) - f(q) + f(q) - f(0)| \geq |f(q) - f(0)| - |f(p) - f(q)| \\ &> (m - \epsilon)|q| - 2\epsilon M|p| = (m - \epsilon)|q - p + p| - 2\epsilon M|p| \\ &\geq (m - \epsilon)(|p| - |q - p|) - 2\epsilon M|p| \geq (m - \epsilon)(|p| - 2\epsilon|p|) - 2\epsilon M|p| \\ &= (m - \epsilon)(1 - 2\epsilon - 2\epsilon M)|p| > \tfrac{m}{2}|p| \end{aligned}$$

for $\epsilon$ sufficiently small. Thus, with such an $\epsilon$ and a corresponding $\delta$, $f(p) \neq f(0)$ for $p \in \Delta(\delta) \cap D$. If $w_0 \in \mathcal{I}(f)$ then there is a sequence $\{z_n\}$ of points in $D$ for $f(z_n) = w_0$. Some subsequence $\{z_{n_i}\}$ converges to some $\zeta' \in \partial D$. Since $f(\zeta') = w_0$ by the continuity of $f$ on $\overline{D}$, what we have just shown implies that $\zeta' \notin \partial D \setminus A$. Thus $\mathcal{I}(f) \subset f(\partial D \setminus A)$. Since $f$ is locally Lipschitz continuous on $\partial D$ and $\lambda_1(\partial D \setminus A) = 0$, it follows that $\lambda_1(\mathcal{I}(f)) = 0$.■

## 3. Proof of Theorem A.

**3.1. Preliminary considerations.** The dimensionality property of $\mathcal{I}(f)$ for the function $f$ we will construct is based on the corresponding property of Cantor-type sets. Let $N \geq 3$ be an integer and let $\alpha \in (0, 1)$. Let

$$E_1 = E_1^{(N, \alpha)} = \cup \{[\tfrac{k}{N} - \tfrac{\alpha}{2N}, \tfrac{k}{N} + \tfrac{\alpha}{2N}] : 1 \leq k \leq N - 1\}.$$

This closed set is the union of the $N - 1$ closed intervals $I_k^{(1)}$ of length $\tfrac{\alpha}{N}$ centered at the points $\tfrac{k}{N}$, $1 \leq k \leq N - 1$, respectively. The next set $E_2$ will be the union of $(N - 1)^2$ closed intervals of length $(\frac{\alpha}{N})^2$ that results by replacing each $I_k^{(1)}$ with $\frac{k}{N} + (\frac{\alpha}{N})(E_1 - \frac{1}{2})$. Continuing in this manner, if we have defined $E_m$, consisting of $(N-1)^m$ disjoint closed intervals $I_k^{(m)}$ of length $(\frac{\alpha}{N})^m$ centered at points $x_1 < \ldots < x_{(N-1)^m}$, then $E_{m+1}$ will be the set consisting of $(N-1)^{m+1}$ closed intervals of length $(\frac{\alpha}{N})^{m+1}$ that results by replacing each of the $I_k^{(m)}$ by $x_k + (\frac{\alpha}{N})^m(E_1 - \frac{1}{2})$. In this manner, we arrive at a nested sequence $\{E_k^{(N,\alpha)}\}$ of closed subsets of $[0,1]$. It can be shown that the $s$-dimensional Hausdorff measure $H^{(s)}(E^{(N,\alpha)})$ of $E^{(N,\alpha)} = \cap \{E_k^{(N,\alpha)} : k \geq 1\}$, satisfies $H^{(s)}(E^{(N,\alpha)}) > 0$, where $s = \frac{\log(N-1)}{\log(\frac{N-1}{\alpha})}$, so that for any given $\alpha \in (0,1)$, $\dim(E^{(N,\alpha)}) \to 1$ as $N \to \infty$. For appropriate sequences $\{\alpha_N\}$, $\{b_N\}$ and $\{\epsilon_N\}$, $\underset{N \geq 3}{\cup} (\epsilon_N E^{(N,\alpha_N)} + b_N)$, $N \geq 1$, will have dimension 1. The set that actually arises in the proof of Theorem A is slightly different from $E^{(N,\alpha)}$, and for the sake of completeness we include a proof of the required dimensionality property, which is simply a reworking of the proof given in [F] that the standard Cantor set has dimension $\frac{\log 2}{\log 3}$. (See also [P, Chapter 10] for similar dimensionality calculations.)

Let $\mathbb{S}$ denote the vertical strip $\mathbb{S} = \{z : 0 < \Re\{z\} < 1\}$. For any measurable $X \subset \mathbb{R}$ for which $\lambda_1(X) < \infty$ we define

(1) $$\mathbf{P}(X, z) = \frac{i}{\pi} \int_X \frac{1}{z-t} \, dt, \, z \in \mathbb{H}.$$

We point out a number of simple properties of the function $\mathbf{P}(X, z)$ which are immediate consequences of the fact that $\mathbf{P}(X, z)$ is to withinin an imaginary constant the Poisson integral of the characteristic function of $X$. First of all, for any such $X \neq \emptyset$, $\mathbf{P}(X, \mathbb{H}) \subset \mathbb{S}$, as long as $\lambda_1(X) \neq 0$. Furthermore, $\mathbf{P}(X, z)$ has a continuous extension to $\mathbb{R} \backslash \partial X$, where $\partial X$ denotes the boundary of $X$ as a subset of $\mathbb{R}$. More specifically, for $x \in \mathbb{R} \backslash \partial X$, $\Re\{\mathbf{P}(X, x)\} = 1$ or 0, respectively, when $x \in X$ or $x \notin X$. Moreover, if $a < b < c$ and $(a, b) \subset X$ and $(b, c) \cap X = \emptyset$, then $\Im\{\mathbf{P}(X, z)\} \to \infty$, $z \to b$ in $\overline{\mathbb{H}}$ and if $(b, c) \subset X$ and $(a, b) \cap X = \emptyset$, then $\Im\{\mathbf{P}(X, z)\} \to -\infty$, $z \to b$ in $\overline{\mathbb{H}}$. If for $z \in \overline{\mathbb{H}} \backslash \{0\}$, we denote by $\log z$ the principal value of the logarithm of $z$, that is, for such $z$, $\Im\{\log z\} \in [0, \pi]$, then for $a < b$, $z \in \mathbb{H} \backslash \{a, b\}$

(2) $$\mathbf{P}((a,b), z) = \frac{i}{\pi} \int_a^b \frac{1}{z-t} \, dt = \frac{i}{\pi}(\log(z-a) - \log(z-b)),$$

from which it follows that $\mathbf{P}((a, b), z)$ maps $\mathbb{H}$ one-to-one onto $\mathbb{S} = \{z : 0 < \Re\{z\} < 1\}$ with

$$\mathbf{P}((a,b), a) = -\infty i, \quad \mathbf{P}((a,b), b) = \infty i, \quad \mathbf{P}((a,b), \infty) = 0 \text{ and } \mathbf{P}((a,b), \tfrac{a+b}{2}) = 1.$$

Under this mapping $(a, b)$ and $\mathbb{R} \backslash [a, b]$ correspond, respectively, to the vertical boundary lines $i\mathbb{R} + 1$ and $i\mathbb{R}$ of $\partial \mathbb{S}$.

Let

$$H_0(z) = \mathbf{P}((-1,1), z) = \tfrac{i}{\pi}(\log(z+1) - \log(z-1)),$$

that is, the one-to-one mapping of $\mathbb{H}$ onto $\mathbb{S}$, with

$$H_0(-1) = -i\infty,\ H_0(1) = i\infty,\ H_0(\infty) = 0 \text{ and } H_0(0) = 1,$$

so that

$$H_0^{-1}(-i\infty) = -1,\ H_0^{-1}(i\infty) = 1,\ H_0^{-1}(0) = \infty \text{ and } H_0^{-1}(1) = 0.$$

A linear change of variable shows that if $a > 0$ and $c$ are real numbers, then

(3) $$\mathbf{P}(aX + c, z) = \mathbf{P}(X, \tfrac{z-c}{a}).$$

With reference to a given function $g$ on $\mathbb{H}$ (and it will always be clear to which such $g$ we are referring), for any $X \subset \mathbb{R}$ of finite measure we define

(4) $$G(X, z) = \int_0^z g'(H_0^{-1}(\mathbf{P}(X, \zeta)))\, d\zeta$$

and use essentially the same notation for the corresponding operation on functions $u$ for which $u(\mathbb{H}) \subset \mathbb{S}$:

(5) $$G(u, z) = \int_0^z g'(H_0^{-1}(u(\zeta)))d\zeta,$$

so that for $X \subset \mathbb{R}$, we have $G(X, z) = G(P_X, z)$, where $P_X(z) = \mathbf{P}(X, z)$.

In §3.2 we shall construct for each $N \geq 3$ a countable family of disjoint open subintervals of $(0, 1)$ for whose union $X^{(N)}$ it will be subsequently shown that if $g$ satisfies the hypothesis of Theorem A, then $\dim(\mathcal{I}(G_N)) \geq d(N) \to 1$ as $N \to 1$, where $G_N(z) = G(X^{(N)}, z)$. The exact expression for $d(N)$ is given in (26) at the end of §3.3.

### 3.2. Construction of $X^{(N)}$.

In what follows $\mathcal{K}$ will denote the set of all finite sequences (to be referred to as *nodes*) $\mathbf{k} = (k_1, k_2, \ldots k_l)$ of integers in $[1, N-1]$. The length of $\mathbf{k}$ will be denoted by $L(\mathbf{k})$ and $\mathcal{K}_l$ will stand for the set of all $\mathbf{k} \in \mathcal{K}$ of length $l$. To facilitate the discussion we introduce an ancestor-descendant relationship: if $\mathbf{k}$ is a proper initial subsequence of $\mathbf{l}$ we say that $\mathbf{k}$ is an *ancestor* of $\mathbf{l}$, and that $\mathbf{l}$ is a *descendant* of $\mathbf{k}$. Nodes $\mathbf{l}_1$ and $\mathbf{l}_2$ are *siblings* if they differ only in their last element; the set of the $N - 2$ siblings of $\mathbf{k}$ will be denoted by $\mathcal{S}(\mathbf{k})$. Node $\mathbf{k}$ is the *parent* of $\mathbf{l}$, if $\mathbf{k}$ is the ancestor of $\mathbf{l}$ for which $L(\mathbf{k}) = L(\mathbf{l}) - 1$. For each $\mathbf{k} \in \mathcal{K}_l$ we define

(6) $$\mathcal{A}(\mathbf{k}) = \{\mathbf{a}_l, \mathbf{a}_{l-1}, \ldots, \mathbf{a}_1\},$$

where $\mathbf{a}_l = \mathbf{k}$ and $\mathbf{a}_i$ is the sole ancestor of $\mathbf{k}$ with $L(\mathbf{a}_i) = i$, that is, $\mathbf{a}_{l-1}$ is the parent of $\mathbf{k}$, $\mathbf{a}_{l-2}$ is its grandparent, etc. Furthermore, we denote by $\mathcal{D}(\mathbf{k})$ the set of descendants of $\mathbf{k}$ and let $\mathcal{C}(\mathbf{k}) = \mathcal{K}_l \backslash \{\mathbf{k}\}$.

We now introduce several positive parameters (in addition to the integer $N \geq 3$) with which we will be working together with certain initial upper bounds on their respective values which subsequent developments in the construction will be replaced with possibly smaller (positive) numbers (these further restrictions will depend solely on the function $g$ of Theorem A). They are $\epsilon < \frac{1}{100}$, $\beta_1 < \frac{1}{100}$ and $\gamma_1 < \epsilon \beta_1$. We impose the additional restriction that $\gamma_1$ be so small that if $U$ is any subset of $(-\frac{1}{2}\beta_1, \frac{1}{2}\beta_1)$ of measure at most $\gamma_1$, then

$$|\mathbf{P}(U, z)| < \epsilon \text{ on } \overline{\mathbb{H} \cap (\Delta \backslash \Delta(\beta_1))}.$$

Furthermore, for notational convenience we abbreviate

$$\alpha = \alpha(\epsilon) = \frac{\epsilon}{N(1+\log N)} \quad \text{and} \quad \gamma = \tfrac{1}{2}\gamma_1 \alpha.$$

and define

$$J = (-1, 1),\ B = (-\beta_1, \beta_1),\ \text{and}\ I = (-1, 1) \backslash \overline{B}.$$

When working with *copies* of $J$ and $I$ of the form $aJ + b$ and $aI + b$, we refer to $aB + b$ as the corresponding *gap*.

We will construct the set $X^{(N)}$, mentioned in the last paragraph of the preceding section, as the union of a set $\{I(\mathbf{k}) : \mathbf{k} \in \mathcal{K}\}$ of copies of $I$. For $\mathbf{k} = (j)$, that is, for a sequence of length 1 whose sole element is the integer $j$, we define

(7) $\quad c(\mathbf{k}) = \frac{j}{N},\ I(\mathbf{k}) = \frac{\alpha}{2}I + c(\mathbf{k}),\ J(\mathbf{k}) = \frac{\alpha}{2}J + c(\mathbf{k})\ \text{and}\ B(\mathbf{k}) = \frac{\alpha}{2}B + c(\mathbf{k})$

and let $X_1 = \bigcup\limits_{L(\mathbf{k})=1}^{N-1} I(\mathbf{k})$. Inside the gap of each of these $I(\mathbf{k})$ we will insert a translate of $\gamma X_1$ centered at a point $x(\mathbf{k})$ near the center $c(\mathbf{k})$ of $I(\mathbf{k})$, then inside the gap of each of the resulting $(N-1)^2$ copies $I(\mathbf{l})$, $l = 2$, of $I$ we will insert a translate of $\gamma^2 X_1$, centered at a point $x(\mathbf{l})$ lying near the center $c(\mathbf{l})$ of $I(\mathbf{l})$, and so on. More precisely, for a node $\mathbf{k}'$ of length greater than 1 consisting of the sequence $\mathbf{k}$ followed by the integer $j$ we define, with $l = L(\mathbf{k})$

(8) $\quad c(\mathbf{k}') = x(\mathbf{k}) + (\frac{j}{N} - \frac{1}{2})\gamma^l,\ I(\mathbf{k}') = \frac{\alpha}{2}\gamma^l I + c(\mathbf{k}'),\ J(\mathbf{k}') = \frac{\alpha}{2}\gamma^l J + c(\mathbf{k}'),$
$\quad\quad B(\mathbf{k}') = \frac{\alpha}{2}\gamma^l B + c(\mathbf{k}').$

The points $x(\mathbf{k}) : \mathbf{k} \in \mathcal{K}$ will be defined below by induction on the length of $\mathbf{k}$ and for the moment we mention that the $x(\mathbf{k})$ will satisfy

(9) $$x(\mathbf{k}) \in \tfrac{\alpha}{8}\gamma^{L(\mathbf{k})-1}B + c(\mathbf{k}), \quad \mathbf{k} \in \mathcal{K},$$

so that, in particular $x(\mathbf{k})$, will lie well inside the gap $B(\mathbf{k})$. It is not the case, however, that any set of points $\{x(\mathbf{k})\}$ satisfing (9) will suffice for our construction since for our purposes it is necessary that there be some positive constant $C_0$ which does not depend either on $N$ or on $\mathbf{k}$ such that

$$|\mathbf{P}(X_{\mathcal{K}\setminus\{\mathbf{k}\}}, z)| < C_0\epsilon, \ z \in c(\mathbf{k}) + \tfrac{\alpha}{2}\gamma^{l-1}(\overline{\mathbb{H} \cap (\Delta\setminus\Delta(\beta_1))}), \ \mathbf{k} \in \mathcal{K}.$$

We will, however, show that the $x(\mathbf{k})$ can be chosen in such a way that this holds with $C_0 = 12$. For any given $\mathcal{U} \subset \mathcal{K}$, we shall write

$$X_\mathcal{U} = \bigcup_{\mathbf{l} \in \mathcal{U}} I(\mathbf{l}).$$

In the following lemmas we assume only (7), (8) and (9).

**Lemma 1**. *Let $L(\mathbf{k}) = l$ and $\mathcal{A}(\mathbf{k})$ be as in (6). Then*

$$\Delta(c(\mathbf{a}_1), \tfrac{\alpha}{4}\beta_1) \supset \Delta(c(\mathbf{a}_2), \gamma\tfrac{\alpha}{4}\beta_1) \supset \ldots \supset \Delta(c(\mathbf{a}_l), \gamma^{l-1}\tfrac{\alpha}{4}\beta_1).$$

**Proof**. Indeed, if $2 \leq i \leq l$ and $z \in \Delta(c(\mathbf{a}_i), \gamma^{i-1}\tfrac{\alpha}{4}\beta_1)$, then $|z - c(\mathbf{a}_i)| < \gamma^{i-1}\tfrac{\alpha}{4}\beta_1$, $|c(\mathbf{a}_i) - x(\mathbf{a}_{i-1})| < \tfrac{1}{2}\gamma^{i-1}$ (by (8)), and $|x(\mathbf{a}_{i-1}) - c(\mathbf{a}_{i-1})| \leq \tfrac{\alpha}{8}\gamma^{i-2}\beta_1$ (by (9)), so that

$$\begin{aligned}|z - c(\mathbf{a}_{i-1})| &\leq |z - c(\mathbf{a}_i)| + |c(\mathbf{a}_i) - x(\mathbf{a}_{i-1})| + |x(\mathbf{a}_{i-1}) - c(\mathbf{a}_{i-1})| \\ &< \gamma^{i-1}\tfrac{\alpha}{4}\beta_1 + \tfrac{1}{2}\gamma^{i-1} + \tfrac{\alpha}{8}\gamma^{i-2}\beta_1 \\ &\leq (\gamma\tfrac{\alpha}{4}\beta_1 + \tfrac{1}{2}\gamma + \tfrac{\alpha}{8}\beta_1)\gamma^{i-2} < \tfrac{\alpha}{4}\gamma^{i-2}\beta_1\end{aligned}$$

since $\gamma\tfrac{\alpha}{4}\beta_1 + \tfrac{1}{2}\gamma < \gamma < \alpha\gamma_1 < \alpha\epsilon\beta_1 < \tfrac{\alpha\beta_1}{8}$, provided only that $\epsilon < \tfrac{1}{8}$. Thus indeed $z \in \Delta(c(\mathbf{a}_i), \gamma^{i-1}\tfrac{\alpha}{4}\beta_1)$ implies that $z \in \Delta(c(\mathbf{a}_{i-1}), \gamma^{i-2}\tfrac{\alpha}{4}\beta_1)$. ∎

**Lemma 2**. *Let $L(\mathbf{k}) = l$ and for each sibling $\mathbf{s} \in \mathcal{S}(\mathbf{k})$ of $\mathbf{k}$ let $Y_\mathbf{s}$ be a measurable subset of $J(\mathbf{s})$, then*

$$\Big|\sum_{\mathbf{s} \in \mathcal{S}(\mathbf{k})} \mathbf{P}'(Y_\mathbf{s}, z)\Big| < \tfrac{8\alpha N^2}{\gamma^{l-1}}, \ z \in \Delta(c(\mathbf{k}), \tfrac{\alpha}{2}\gamma^{l-1}).$$

**Proof**. If the last element of the sequence $\mathbf{k}$ is $k$, then for each integer $t \in [1, N-2]$ there are at most two integers between 1 and $N-1$ whose distance from $k$ is $t$, so that there are at most two siblings $\mathbf{s}$ of $\mathbf{k}$ for which the distance between $c(\mathbf{k})$ and $c(\mathbf{s})$ is $t\tfrac{\gamma^{l-1}}{N}$, as follows from (8) of the $c(\mathbf{l})$. For such $\mathbf{s}$ the minimum distance between points of $\Delta(c(\mathbf{k}), \tfrac{\alpha}{2}\gamma^{l-1})$ and $J(\mathbf{s})$ is $(\tfrac{t}{N} - \alpha)\gamma^{l-1} \geq (t-\epsilon)\tfrac{\gamma^{l-1}}{N} \geq \tfrac{t}{2}\tfrac{\gamma^{l-1}}{N}$, since $\epsilon < 1$. From this and the fact that $\lambda_1(J(\mathbf{s})) = \alpha\gamma^{l-1}$ it follows that for $z \in \Delta(c(\mathbf{k}), \alpha\gamma^{l-1})$

$$|\mathbf{P}'(Y_{\mathbf{s}}, z)| = \tfrac{1}{\pi}|\int_{J(\mathbf{s})} \tfrac{dt}{(z-t)^2}| \leq \tfrac{1}{\pi}\int_{J(\mathbf{s})} \tfrac{dt}{|(z-t)^2|} \leq \tfrac{\alpha}{\pi}\gamma^{l-1}\tfrac{(2N)^2}{(t\gamma^{l-1})^2} = \tfrac{4\alpha}{\pi}\tfrac{N^2}{t^2\gamma^{l-1}},$$

for $z \in \Delta(c(\mathbf{k}), \tfrac{\alpha}{2}\gamma^{l-1})$. Thus,

$$|\sum_{\mathbf{s}\in\mathcal{S}(\mathbf{k})} \mathbf{P}'(Y_{\mathbf{s}}, z)| \leq \tfrac{2}{\pi}\tfrac{4\alpha N^2}{\gamma^{l-1}}\sum_{t=1}^{N-2}\tfrac{1}{t^2} < \tfrac{4\alpha N^2}{\gamma^{i-1}}\sum_{t=1}^{\infty}\tfrac{1}{t^2} < \tfrac{8\alpha N^2}{\gamma^{i-1}},$$

as claimed since $\sum_{t=1}^{\infty}\tfrac{1}{t^2} = \tfrac{\pi^2}{6} < 2$. ∎

**Lemma 3**. *Let $L(\mathbf{k}) = l$ and $\mathcal{A}(\mathbf{k})$ be as in (6). Then,*

$$|\mathbf{P}'(X_{\mathcal{A}(\mathbf{k})}, z)| < \tfrac{4}{\alpha\beta_1\gamma^{l-1}}, \; z \in \Delta(c(\mathbf{k}), \tfrac{\alpha\beta_1}{4}\gamma^{l-1}).$$

**Proof**. Clearly, from (2) we have

(10) $\quad \mathbf{P}'(I, z) = \tfrac{i}{\pi}\{\tfrac{1}{z+1} - \tfrac{1}{z-1} - \tfrac{1}{z+\beta_1} + \tfrac{1}{z-\beta_1}\} = \tfrac{2i}{\pi}\tfrac{(\beta_1-1)(z^2+\beta_1)}{(1-z^2)(\beta_1^2-z^2)},$

so that, since $\beta_1 < \tfrac{1}{100}$, we have for $|z| \leq \tfrac{1}{2}\beta_1$ that

$$|\mathbf{P}'(I, z)| \leq \tfrac{2}{\pi}\big(\tfrac{(1-\beta_1)(\beta_1+\tfrac{1}{4}\beta_1^2)}{(1-\tfrac{1}{4}\beta_1^2)(\tfrac{3}{4}\beta_1^2)}\big) < \tfrac{8}{3\pi}\tfrac{1+\tfrac{1}{4}\beta_1}{(1-\tfrac{1}{4}\beta_1^2)\beta_1} \leq \tfrac{1}{\beta_1}.$$

Since by (3) of §3.1, $\mathbf{P}(I(\mathbf{a}_i), z) = \mathbf{P}(\tfrac{\alpha}{2}\gamma^{l-1}I + c(\mathbf{a}_i), z) = \mathbf{P}(I, \tfrac{z-c(\mathbf{a}_i)}{\tfrac{\alpha}{2}\gamma^{l-1}})$, it follows that

$$|\mathbf{P}'(I(\mathbf{a}_i), z)| \leq \tfrac{2}{\alpha\beta_1\gamma^{i-1}}, \; z \in \Delta(c(\mathbf{a}_i), \tfrac{\alpha\beta_1}{4}\gamma^{i-1}), 1 \leq i \leq l.$$

From this together with nesting of the disks $\Delta(c(\mathbf{a}_i), \tfrac{\alpha\beta_1}{4}\gamma^{i-1})$ established in Lemma 1, and the fact that $\gamma < \tfrac{1}{2}$ we conclude that

$$|\mathbf{P}'(X_{\mathcal{A}(\mathbf{k})}, z)| = |\sum_{i=1}^{l} \mathbf{P}'(I(\mathbf{a}_i), z)| \leq \sum_{i=1}^{l}\tfrac{2}{\alpha\beta_1\gamma^{i-1}} < \tfrac{2}{\alpha\beta_1\gamma^{l-1}}\sum_{i=0}^{\infty}\gamma^i < \tfrac{4}{\alpha\beta_1\gamma^{l-1}},$$

for $z \in \Delta(c(\mathbf{k}), \tfrac{\alpha\beta_1}{4}\gamma^{l-1})$ as claimed. ∎

**Lemma 4**. *Let $L(\mathbf{k}) = l$. For each node $\mathbf{m} \in \mathcal{C}(\mathbf{k}) = \mathcal{K}_l \backslash \{\mathbf{k}\}$, let $Y_{\mathbf{m}}$ be a measurable subset of $J(\mathbf{m})$, then*

$$|\sum_{\mathbf{m}\in\mathcal{C}(\mathbf{k})} \mathbf{P}(Y_{\mathbf{m}}, z)| < 4\epsilon, \; z \in \overline{\mathbb{H} \cap \Delta(c(\mathbf{k}), \tfrac{\alpha}{2}\gamma^{l-1})}.$$

**Proof.** As in (6) let $\mathcal{A}(\mathbf{k}) = \{\mathbf{a}_l = \mathbf{k}, \mathbf{a}_{l-1}, \ldots, \mathbf{a}_1\}$. For $1 \leq i \leq l-1$, $\mathcal{L}_i$ be the set of descendants of the siblings of $\mathbf{a}_i$ of length $l$ (that is, $\mathcal{L}_i = \mathcal{D}(\mathbf{a}_i) \cap \mathcal{K}_l$) and let $\mathcal{L}_l = \mathcal{S}(\mathbf{k})$, the set of siblings of $\mathbf{k}$. Clearly,

$$\mathcal{C}(\mathbf{k}) = \bigcup_{i=1}^{l} \mathcal{L}_i,$$

so that

$$\left| \sum_{\mathbf{m} \in \mathcal{C}(\mathbf{k})} \mathbf{P}(Y_{\mathbf{m}}, z) \right| = \left| \sum_{i=1}^{l} \sum_{\mathbf{m} \in \mathcal{L}_i} \mathbf{P}(Y_{\mathbf{m}}, z) \right|.$$

We estimate each of the inner sums $\sigma_i = \sum_{\mathbf{m} \in \mathcal{L}_i} \mathbf{P}(Y_{\mathbf{m}}, z)$. For each of the $N-2$ siblings $\mathbf{b}_1, \ldots, \mathbf{b}_{N-2}$ of $\mathbf{a}_i$ there are $(N-1)^{l-i} < N^{l-i}$ summands in $\sigma_i$. Let the last element of the sequence $\mathbf{a}_i$ be $k$. As pointed out in the proof of Lemma 2, for each integer $t \in [1, N-2]$ there are at most 2 integers $j$ between 1 and $N-1$, at distance $t$ from $k$, so that for such $j$ the distance from $c(\mathbf{b}_j)$ and $c(\mathbf{a}_i)$ is $t\frac{\gamma^{i-1}}{N}$ since for nodes of length $i$ the distances between the $c$'s corresponding to siblings are multiples of $\frac{\gamma^{i-1}}{N}$. For any descendant $\mathbf{m} \in \mathcal{C}(\mathbf{k})$ of such a $\mathbf{b}_j$ the minimum distance between points of $\Delta(c(\mathbf{a}_i), \frac{\alpha}{2}\gamma^{i-1})$ and $J(\mathbf{m})$ is $(t-\alpha)\frac{\gamma^{i-1}}{N} \geq (\frac{t}{N} - \epsilon)\gamma^{i-1} \geq \frac{t}{2}\frac{\gamma^{i-1}}{N}$, since $\epsilon < \frac{1}{2}$. Grouping the $(N-1)^{l-i}$ summands $\sigma_i$ according as they are descendants of $\mathbf{b}_1, \ldots, \mathbf{b}_{N-2}$ and taking into account that the linear measure of $Y_{\mathbf{m}}$ for each of these descendants $\mathbf{m}$ is at most $\alpha\gamma^{l-1}$

$$\left| \sum_{\mathbf{m} \in \mathcal{L}_i} \mathbf{P}(Y_{\mathbf{m}}, z) \right| \leq \frac{2}{\pi} N^{l-i} \alpha \gamma^{l-1} \sum_{t=1}^{N-2} \frac{2N}{t\gamma^{i-1}}$$

$$= \frac{4}{\pi} \alpha N^{l+1-i} \gamma^{l-i} \sum_{t=1}^{N-2} \frac{1}{t} < 2\alpha N^{l-i} \gamma^{l-i} N \log N = 2\alpha (N\gamma)^{l-i} N \log N,$$

for $z \in \overline{\mathbb{H} \cap \Delta(c(\mathbf{a}_i), \frac{\alpha}{2}\gamma^{i-1})}$. But $\alpha N \log N < \epsilon$, and $N\gamma = \frac{1}{2}N\alpha\gamma_1 < \frac{1}{2}$ and therefore,

$$\left| \sum_{\mathbf{m} \in \mathcal{C}(\mathbf{k})} \mathbf{P}(Y_{\mathbf{m}}, z) \right| < 2\epsilon \sum_{i=1}^{l} \frac{1}{2^{l-i}} < 4\epsilon,$$

as asserted. ∎

**Lemma 5.** *Let $L(\mathbf{k}) = l$ and let $Y \subset c(\mathbf{k}) + \frac{\alpha}{4}\gamma^{l-1}B$ have measure at most $\gamma^l$. Then $|\mathbf{P}(Y, z)| < \epsilon$ for $z \in c(\mathbf{k}) + \frac{\alpha}{2}\gamma^{l-1}(\overline{\mathbb{H} \cap (\Delta \backslash \Delta(\beta_1))})$.*

**Proof.** For convenience we abbreviate $\frac{\alpha}{2}\gamma^{l-1}$ by $\alpha'$ and $c(\mathbf{k})$ by $c$. Let $W = \frac{Y-c}{\alpha'}$. From (3) in §3.1, $\mathbf{P}(Y, z) = \mathbf{P}(W, \frac{z-c}{\alpha'})$. But $z \in c(\mathbf{k}) + \frac{\alpha}{4}\gamma^{l-1}B$ is equivalent to $\frac{z-c}{\alpha'} \in \frac{1}{2}B$, so that $W \subset \frac{1}{2}B$. Also, $\lambda_1(W) = \frac{1}{\alpha'}\lambda_1(Y) < \gamma^l / \frac{\alpha}{2}\gamma^{l-1} = \frac{2}{\alpha}\gamma = \gamma_1$, so that from the definition of $\gamma_1$, $|\mathbf{P}(W, \frac{z-c}{\alpha'})| < \epsilon$ for $\frac{z-c}{\alpha'} \in \overline{\mathbb{H} \cap (\Delta \backslash \Delta(\beta_1))}$, i.e., for $z \in c(\mathbf{k}) + \frac{\alpha}{2}\gamma^{l-1}(\overline{\mathbb{H} \cap (\Delta \backslash \Delta(\beta_1))})$. ∎

For convenience we define

$$\mathcal{K}'(\mathbf{k}) = \{\mathbf{l} : L(\mathbf{l}) \leq L(\mathbf{k})\}.$$

We will now show by a simple induction on $L(\mathbf{k})$ that we can choose the points $x(\mathbf{k})$, $\mathbf{k} \in \mathcal{K}$ in such a way that when the $c(\mathbf{k}), J(\mathbf{k}), I(\mathbf{k}), B(\mathbf{k})$ aredefined by (7) and (8), then we have

(11) $$\mathbf{P}(X_{\mathcal{K}'(\mathbf{k})}, x(\mathbf{k})) = 0,$$

(12) $$x(\mathbf{k}) \in 4\epsilon\alpha\gamma^{L(\mathbf{k})-1}B + c(\mathbf{k}), \text{ i.e., } |x(\mathbf{k}) - c(\mathbf{k})| < 4\epsilon\alpha\beta_1\gamma^{L(\mathbf{k})-1}$$

(13) $$|\mathbf{P}(X_{\mathcal{K}'(\mathbf{k})\backslash\{\mathbf{k}\}}, z)| < 7\epsilon, \ z \in \overline{\mathbb{H} \cap \Delta(c(\mathbf{k}), \frac{\alpha}{2}\gamma^{L(\mathbf{k})-1})}.$$

To begin the induction let $l = 1$, so that $\mathbf{k} = (j)$, for some $j \in [1, N-1]$. It follows from (7) and Lemma 4 with $l = 1$ and $Y_{\mathbf{m}} = I(\mathbf{m}) \subset J(\mathbf{m})$, for each $\mathbf{m} \in \mathcal{K}'(\mathbf{k}) \backslash \{\mathbf{k}\}$, that

(14) $$\mathbf{P}(X_{\mathcal{K}'(\mathbf{k})\backslash\{\mathbf{k}\}}, z) < 4\epsilon, \ z \in \overline{\mathbb{H} \cap \Delta(c(\mathbf{k}), \frac{\alpha}{2})}.$$

Now, as indicated in the comments immediately following the definition (1) of $\mathbf{P}(X, z)$, $\mathbf{P}(X_{\mathcal{K}_1}, z)$ is pure imaginary on each of the gaps $B(\mathbf{k}) = \frac{\alpha}{2}B + c(\mathbf{k}) = (\frac{j}{N} - \frac{\alpha}{2}\beta_1, \frac{j}{N} + \frac{\alpha}{2}\beta_1)$. It follows from the formula (10) for $\mathbf{P}'(I, z)$ in the proof of Lemma 3 that

(15) $$\Im\{\mathbf{P}'(I, z)\} < -\frac{1}{2\beta_1}, z \in B.$$

Since $\mathbf{P}(I(\mathbf{k}), z) = \mathbf{P}(\frac{\alpha}{2}I + c(\mathbf{k}), z) = \mathbf{P}(I, \frac{z-c(\mathbf{k})}{\alpha/2})$, $\Im\{\mathbf{P}'(I(\mathbf{k}), z)\} < -\frac{1}{\alpha\beta_1}$ on $B(\mathbf{k})$. Since $\mathbf{P}(I(\mathbf{k}), c(\mathbf{k})) = 0$, it follows from (14) and the intermediate value theorem that there indeed exists an $x(\mathbf{k})$ within $4\epsilon\alpha\beta_1$ of $c(\mathbf{k})$ for which (11) holds, so that condition (12) is satisfied. Finally, (13) follows immediately from (14).

Now assume that $l \geq 1$ and that we have defined $x(\mathbf{k})$ for all $\mathbf{k}$ of length at most $l$ in such a way that (11)-(13) hold. Let $L(\mathbf{k}) = l$, let $j$ be an integer in $[1, N-1]$ and let $\mathbf{k}'$ denote the sequence of length $l+1$ consisting of $\mathbf{k}$ followed by $j$. The objects $c(\mathbf{k}')$, $I(\mathbf{k}'), J(\mathbf{k}')$ and $B(\mathbf{k}')$ are as defined in (8). Let $\mathcal{A}(\mathbf{k}) = \{\mathbf{a}_l = \mathbf{k}, \mathbf{a}_{l-1}, \ldots, \mathbf{a}_1\}$, the sequence of $\mathbf{k}$ together with its ancestors, as in (6). It is clear that $\Delta(c(\mathbf{a}_i), \frac{\alpha}{2}\gamma^{i-1}) \supset \Delta(c(\mathbf{a}_{i+1}), \frac{\alpha}{2}\gamma^i)$, $1 \leq i \leq l-1$ so that from Lemma 2 it follows that

$$|\mathbf{P}'(X_{\mathcal{K}'(\mathbf{k})\backslash\mathcal{A}(\mathbf{k})}, z)| < 8\alpha N^2 \sum_{i=1}^{l} \frac{1}{\gamma^{i-1}}, z \in \Delta(c(\mathbf{k}), \frac{\alpha}{2}\gamma^{l-1}),$$

so that,

(16) $\quad |\mathbf{P}(X_{\mathcal{K}'(\mathbf{k})\setminus\mathcal{A}(\mathbf{k})}, z_1) - \mathbf{P}(X_{\mathcal{K}'(\mathbf{k})\setminus\mathcal{A}(\mathbf{k})}, z_2)| < 8\alpha N^2 |z_1 - z_2| \sum_{i=1}^{l} \frac{1}{\gamma^{i-1}},$

for $z_1, z_2 \in \overline{\mathbb{H} \cap \Delta(c(\mathbf{k}), \frac{\alpha}{2}\gamma^{l-1})}$. But $|x(\mathbf{k}) - c(\mathbf{k})| < 4\epsilon\alpha\beta_1\gamma^{l-1}$, so that

(17) $\quad \Delta(x(\mathbf{k}), \gamma^l) \subset \Delta(c(\mathbf{k}), 4\epsilon\alpha\beta_1\gamma^{l-1} + \gamma^l) \subset \Delta(c(\mathbf{k}), \frac{\alpha\beta_1}{4}\gamma^{l-1}),$

since $4\epsilon\alpha\beta_1\gamma^{l-1} + \gamma^l = (4\epsilon\alpha\beta_1 + \gamma)\gamma^{l-1} = (4\epsilon\alpha\beta_1 + \gamma)\gamma^{l-1} = (4\epsilon\alpha\beta_1 + \frac{\alpha\gamma_1}{2})\gamma^{l-1}$
$= \alpha(4\epsilon\beta_1 + \frac{\epsilon\beta_1}{2})\gamma^{l-1} = \frac{9}{2}\epsilon\alpha\beta_1\gamma^{l-1} < \frac{\alpha\beta_1}{4}\gamma^{l-1}$. Thus, (16) is valid for $z_1, z_2 \in \overline{\mathbb{H} \cap \Delta(x(\mathbf{k}), \gamma^l)}$, since $\frac{\alpha\beta_1}{4} < \frac{\alpha}{2}$. By Lemma 3 we have

$$|\mathbf{P}'(X_{\mathcal{A}(\mathbf{k})}, z)| < \frac{4}{\alpha\beta_1\gamma^{l-1}}, \ z \in \Delta(c(\mathbf{k}), \frac{\alpha\beta_1}{4}\gamma^{l-1}),$$

so that

(18) $\quad |\mathbf{P}(X_{\mathcal{A}(\mathbf{k})}, z_1) - \mathbf{P}(X_{\mathcal{A}(\mathbf{k})}, z_2)| < \frac{4|z_1 - z_2|}{\alpha\beta_1\gamma^{l-1}},$

for $z_1, z_2 \in \overline{\mathbb{H} \cap \Delta(c(\mathbf{k}), \frac{\alpha\beta_1}{4}\gamma^{l-1})}$. In light of (17), (18) is valid for $z_1$, $z_2 \in \overline{\mathbb{H} \cap \Delta(x(\mathbf{k}), \gamma^l)}$. From (16) and (18) it therefore follows that for $z \in \overline{\mathbb{H} \cap \Delta(x(\mathbf{k}), \gamma^l)}$

$$|\mathbf{P}(X_{\mathcal{K}'(\mathbf{k})}, z) - \mathbf{P}(X_{\mathcal{K}'(\mathbf{k})}, x(\mathbf{k}))| < \frac{4\gamma}{\alpha\beta_1} + 8\alpha N^2\gamma^l \sum_{i=1}^{l} \frac{1}{\gamma^{i-1}}.$$

But since $\gamma = \frac{1}{2}\gamma_1\alpha < \frac{1}{2}$ and $\gamma < \epsilon\beta_1$,

$$8\alpha N^2\gamma^l \sum_{i=1}^{l} \frac{1}{\gamma^{i-1}} < 8\alpha N^2\gamma \sum_{i=0}^{\infty} \gamma^i < 16\alpha N^2\gamma = 8(\alpha N)^2\gamma_1$$
$$= 8\left(\frac{\epsilon}{(1+\log N)}\right)^2 \gamma_1 < 8\left(\frac{\epsilon}{(1+\log N)}\right)^2 \epsilon\beta_1 < \epsilon.$$

Also $\frac{4\gamma}{\alpha\beta_1} = 2\frac{\gamma_1}{\beta_1} < 2\epsilon$. Thus by (11),

(19) $\quad |\mathbf{P}(X_{\mathcal{K}'(\mathbf{k})}, z)| < 3\epsilon, \ z \in \overline{\mathbb{H} \cap \Delta(x(\mathbf{k}), \gamma^l)}.$

All of the $c(\mathbf{l}), J(\mathbf{l}), I(\mathbf{l})$ with $L(\mathbf{l}) = l + 1$ are defined by (8) so that we can address (13) before defining $x(\mathbf{l})$ and verifying (11) and (12). By Lemma 4

(20) $\quad |\sum_{\mathbf{l} \in \mathcal{C}(\mathbf{k}')} \mathbf{P}(I(\mathbf{l}), z)| < 4\epsilon, \ z \in \overline{\mathbb{H} \cap \Delta(c(\mathbf{k}'), \frac{\alpha}{2}\gamma^l)},$

so that, since $\Delta(c(\mathbf{k}'), \frac{\alpha\gamma^l}{2}) \subset \Delta(x(\mathbf{k}), \gamma^l)$, from (19) and (20) we conclude that

(21) $\quad |\sum_{\mathbf{l} \in \mathcal{K}'(\mathbf{k}') \setminus \{\mathbf{k}'\}} \mathbf{P}(I(\mathbf{l}), z)| < 7\epsilon, \ z \in \overline{\mathbb{H} \cap \Delta(c(\mathbf{k}'), \frac{\alpha}{2}\gamma^l)},$

which establishes (13) for $l + 1$. We now show that $x(\mathbf{l})$ can be defined and in a way that makes (11) and (12) valid. This, however, follows exactly as in the case $l = 1$. Since $\mathbf{P}(I(\mathbf{k}'), z) = \mathbf{P}(\frac{\alpha}{2}\gamma^l I + c(\mathbf{k}), z) = \mathbf{P}(I, \frac{z-c(\mathbf{k})}{\alpha\gamma^l/2})$, $\Im\{\mathbf{P}'(I(\mathbf{k}'), z)\} < -\frac{1}{\alpha\gamma^l \beta_1}$ on $B(\mathbf{k}')$ by (15), it follows from (21) and the intermediate value theorem that there is an $x(\mathbf{k}')$ within $7\epsilon\alpha\beta_1\gamma^l$ of $c(\mathbf{k}')$ for which $\mathbf{P}(X_{\mathcal{K}'(\mathbf{k}')}, x(\mathbf{k}')) = 0$, which establishes (11) and (12). Hence by induction we have defined $\{x(\mathbf{k}) : \mathbf{k} \in \mathcal{K}\}$ in such a way that (11), (12) and (13) hold for all $\mathbf{k} \in \mathcal{K}$.

It is clear from Lemma 4 that if

$$\mathcal{N}(\mathbf{k}) = \{\mathbf{l} : L(\mathbf{l}) > L(\mathbf{k}), \mathbf{l} \notin \mathcal{D}(\mathbf{k})\},$$

then for all $\mathbf{k} \in \mathcal{K}$

(22) $\qquad |\mathbf{P}(X_{\mathcal{N}(\mathbf{k})}, z)| < 4\epsilon, \quad z \in \overline{\mathbb{H} \cap \Delta(c(\mathbf{k}), \frac{\alpha}{2}\gamma^{l-1})},$

and follows from Lemma 5 that

$$|\mathbf{P}(X_{\mathcal{D}(\mathbf{k})}, z)| < \epsilon, \ z \in c(\mathbf{k}) + \overline{\mathbb{H} \cap \frac{\alpha}{2}\gamma^{l-1}(\Delta \setminus \Delta(\beta_1))}$$

Putting this together with (13) and (22) we conclude that for

(23) $\qquad |\mathbf{P}(X_{\mathcal{K} \setminus \{\mathbf{k}\}}, z)| < 12\epsilon, z \in c(\mathbf{k}) + \frac{\alpha}{2}\gamma^{l-1}(\overline{\mathbb{H} \cap (\Delta \setminus \Delta(\beta_1))}), \ \mathbf{k} \in \mathcal{K}.$

We define the fundamental sets

$$X^{(N)} = X_{\mathcal{K}} = \bigcup_{\mathbf{k} \in \mathcal{K}} I(\mathbf{k})$$

$$E^{(N)} = \bigcap_{n=1}^{\infty} (\bigcup_{L(\mathbf{l})=n} \overline{J(\mathbf{l})}).$$

### 3.3. Dimension of $E^{(N)}$.

The set $E = E^{(N)}$ is very similar to the set described §3.1 except that we are using the $x(\mathbf{k})$ instead of the $c(\mathbf{k})$, so that things are slightly "off center," but this does not affect the calculation. Let $\{U_i\}$ be a covering of $E$. We can assume that all these $U_i$ are intervals since we are only interested in $\sum_i |U_i|^s$, where here $|U|$ denotes the diameter of $U$. Furthermore, since we can enlarge each $U_i$ in such a way that $\sum_i |U_i|^s$ is increased by an arbitrarily small amount, we can assume that all of the $U_i$ are open and therefore, since $E$ is compact, we can also assume that $\{U_i\}$ is a finite covering. In light of the definition of $s$-dimensional Hausdorff measure we may also assume that $|U_i| < \frac{1}{N} - \alpha$. For each $i$ let $n(i) \geq 1$ be the integer $n$ for which

(24) $$\gamma^{n(i)}(\tfrac{1}{N} - \alpha) \leq |U_i| < \gamma^{n(i)-1}(\tfrac{1}{N} - \alpha).$$

The minimum distance between the intervals which make up $E_l = \bigcup_{L(\mathbf{l})=l} \overline{J(\mathbf{l})}$ is $\gamma^{l-1}(\tfrac{1}{N} - \alpha)$, so that $U_i$ can intersect at most one of the $\overline{J(\mathbf{l})}$ that make up $E_{n(i)}$, and therefore, more generally, for $m \geq 0$, $U_i$ can intersect at most $(N-1)^m$ of the closed intervals which make up $E_{n(i)+m}$ since $E_{n(i)+m} \subset E_{n(i)}$. In other words, if $j \geq n(i)$, then $U_i$ intersects at most $(N-1)^{j-n(i)}$ of the intervals in $E_j$. Now fix $j$ so large that $\gamma^{j+1}(\tfrac{1}{N} - \alpha) \leq |U_i|$ for all $i$. Since $\{U_i\}$ is a covering of $E$, for each interval $\overline{J(\mathbf{l})}$ making up $E_j$, some $U_i$ must touch $\overline{J(\mathbf{l})}$. Thus,

(25) $$(N-1)^j = \#\text{of intervals in } E_j \leq \sum_i (N-1)^{j-n(i)}.$$

We express $N-1$ as $(\tfrac{1}{\gamma})^s$, where $s = \tfrac{\log(N-1)}{\log(1/\gamma)} > 0$. From (24) we have that

$$(N-1)^{-n(i)} = (\gamma^{n(i)})^s \leq (\tfrac{N}{1-N\alpha})^s |U_i|^s,$$

so that from (25) it follows that

$$(N-1)^j \leq \sum_i (N-1)^{j-n(i)} = (N-1)^j (\tfrac{N}{1-N\alpha})^s \sum_i |U_i|^s,$$

and therefore that $\sum_i |U_i|^s \geq (\tfrac{1-N\alpha}{N})^s > 0$. Thus $H^{(s)}(E) > 0$, which means that $\dim(E) \geq s = \tfrac{\log(N-1)}{\log(1/\gamma)}$. Since $\tfrac{1}{\gamma} = \tfrac{2N(1+\log N)}{\epsilon \gamma_1}$, we have that

(26) $$\dim(E^{(N)}) \geq d(N) = \frac{\log(N-1)}{\log N + \log(1+\log N) + \log(2/\epsilon\gamma_1)},$$

which tends to 1 as $N \to \infty$.

### 3.4. A Normalization.

The following lemma allows us to introduce some simplifying conditions.

**Lemma 6**. *Let $g_0$ be nonunivalent in $\mathbb{H}$ with nonvanishing derivative there. Then there is a function $\rho_0$ which maps $\mathbb{H}$ into itself such that $g(z) = \int_0^z g_0'(\rho(\zeta))\,d\zeta$ satisfies*
(G*i*) *$g'$ is analytic and nonvanishing on $\overline{\mathbb{H}}$,*
(G*ii*) *$g'(z)$ is analytic at $\infty$ and $g'(\infty) = b_0 \neq 0$,*
(G*iii*) *$g(z_0) = g(0) = 0$ where $z_0 \in \mathbb{H} \cap \partial\Delta(\tfrac{1}{2})$.*

**Proof**. For $r > 1$, let $q_r$ map $\mathbb{H}$ one-to-one onto $\Delta(ri, r - \tfrac{1}{r})$ with $q_r(i) = i$ and $q_r'(i) > 0$, so that $q_r(z) \to z$ uniformly on compact sets of $\mathbb{H}$ as $r \to \infty$. Thus, for all sufficiently large $r$ the $\int_0^z g'(q_r(\zeta))\,d\zeta$ is analytic on $\overline{\mathbb{H}}$ and nonunivalent on $\mathbb{H}$. We choose such a value of $r$ and denote the corresponding antiderivative by $f(z)$. Since $\overline{q_r(\mathbb{H})}$ is compact, so is $\overline{f'(\mathbb{H})}$, and $0 \notin \overline{f'(\mathbb{H})}$. Also, note that $f'$ is analytic and nonvanishing on $\overline{\mathbb{H}}$ and $f'(1/z)$ is analytic with value $b_0 \neq 0$ at $0$. Let $M'$ be an upper bound for $f'$ on $\mathbb{H}$ and let $|f'(z) - b_0| < \frac{|b_0|}{2}$ for $|z| > R$. Let $z_1, z_2 \in \mathbb{H}$ be distinct with $|z_2| \geq TR$, where $T > 1$. If $[z_1, z_2] \cap \Delta(R) = \emptyset$, then

$$|f(z_2) - f(z_1)| = |\textstyle\int_{[z_1,z_2]}(b_0 + f'(z) - b_0)\,dz| \geq \tfrac{|b_0|}{2}|z_2 - z_1|.$$

If $(z_1, z_2) \cap \Delta(R) = (z_1', z_2') \neq \emptyset$ with $z_1'$ coming before $z_2'$ as we move from $z_1$ to $z_2$ along the segment, then clearly $|z_2 - z_1| \geq (T-1)R$ and

$$\begin{aligned}|f(z_2) - f(z_1)| &= |\textstyle\int_{[z_1,z_2]\setminus[z_1',z_2']}(f'(z) - b_0) + b_0)\,dz + \int_{[z_1',z_2']} f'(z)\,dz| \\ &\geq \tfrac{|b_0|}{2}(|z_2-z_1| - |z_2'-z_1'|) - M'|z_2'-z_1'| \\ &\geq \tfrac{|b_0|}{2}((T-1)R - 2R) - 2M'R = (\tfrac{|b_0|}{2}(T-3) - 2M')R > 0\end{aligned}$$

as long as $T > \frac{4M'}{|b_0|} + 3$. Thus we can only have $f(z_1) = f(z_2)$ for $z_1 \neq z_2$ in $\mathbb{H}$ if $z_1, z_2$ are both in $\Delta(r_0)$, where $r_0 = (\frac{4M'}{|b_0|} + 3)R$. From this it follows that for all sufficiently large $s$, $f$ is univalent on $\overline{\mathbb{H}} + si$. Let $s_0 = \inf\{s : f \text{ is univalent on } \overline{\mathbb{H}} + si\} > 0$. For each $s \in (0, s_0)$ there are distinct points $a_s, b_s \in \Delta(r_0) \cap (\overline{\mathbb{H}} + si)$ for which $f(a_s) = f(b_s)$. From a straightforward compactness argument and the fact that $f$ is locally univalent on $\overline{\mathbb{H}}$ it follows that there exist distinct $a, b \in \mathbb{R} + s_0 i$ whose images under $f$ coincide. Since $f$ is univalent on $\mathbb{H} + s_0 i$, the image of the curve $\mathbb{R} + s_0 i$ is tangent to itself at the point $f(a) = f(b)$ and $f'(b)$ is a negative multiple of $f'(a)$. Therefore, for all sufficiently small $\delta \in (0, s_0)$, with $a' = a - \delta i \in \mathbb{R} + (s_0 - \delta)i$ there holds

$$f(a') \in f((\mathbb{H} + s_0 i)) \subset f(\mathbb{H} + (s_0 - \delta)i)),$$

so that there is a point $b' \in \mathbb{H} + (s_0 - \delta)i$ such that $f(a') = f(b')$. Clearly, $b' - a' \in \mathbb{H}$. Let

$$g(z) = \tfrac{1}{2|b'-a'|}(f(2|b'-a'|z + a') - f(a')),$$

so that $g(0) = g(\frac{b'-a'}{2|b'-a'|}) = 0$. Also, $g'(z)$ is clearly of the form $g_0'(\rho_0(z))$ with $\rho_0(\mathbb{H}) \subset \mathbb{H}$. ∎

Henceforth we shall work with a function $g$ which has the properties given in the conclusion of this lemma.

### 3.5. $G(X^{(N)}, z)$ is bi-Lipschitz on $E^{(N)}$.

We begin with the following

**Lemma 7.** Let $b$, $\eta$, $T > 0$, let $N \geq 3$ be an integer and let $X \subset [0, b]$ satisfy $\lambda_1(X \cap W) \leq \frac{\eta b}{N(1+\log N)}$ for all intervals $W$ of length $\frac{b}{N}$. Then
$$|\mathbf{P}(X, x + \tfrac{Tb\eta}{N}i)| \leq \tfrac{1}{\pi}(\tfrac{3}{T} + 2\eta), \ 0 \leq x \leq b.$$

**Proof.** Let $x \in [0, b]$. For $|k| \leq N$, let $J_k$ be the interval $x + \frac{b}{N}[k - \frac{1}{2}, k + \frac{1}{2}]$. The union of these intervals of length $\frac{b}{N}$ covers $[0, b]$. Let $z = x + \frac{Tb\eta}{N}i$. Since $\frac{1}{|z-t|} \leq \frac{N}{Tb\eta}$, for $t \in \mathbb{R}$, it follows from the hypothesis that
$$\sum_{|k| \leq 1} |\mathbf{P}(X \cap J_k, z)| \leq \tfrac{3\eta b}{\pi N(1+\log N)} \cdot \tfrac{N}{Tb\eta} < \tfrac{3}{\pi T}.$$
.
For $|k| \geq 2$, $\operatorname{dist}(z, J_k) = \frac{b}{N}(|k| - \frac{1}{2}) > \frac{b(|k|-1)}{N}$, so that for $t \in J_k$, $\frac{1}{|z-t|} \leq \frac{N}{b(|k|-1)}$. Thus, it follows from the hypothesis that
$$\sum_{2 \leq |k| \leq N} |\mathbf{P}(X \cap J_k, z)| \leq \tfrac{\eta b}{\pi N(1+\log N)} \sum_{2 \leq |k| \leq N} \tfrac{N}{b(|k|-1)} = \tfrac{2\eta}{\pi(1+\log N)} \sum_{k=1}^{N-1} \tfrac{1}{k} < \tfrac{2\eta}{\pi}.$$

The desired conclusion follows immediately from these two bounds. ∎

We next use Lemma 7 to establish that the mapping $G(X^{(N)}, z)$ is bi-Lipschitz on $E^{(N)}$ and, as will be significant when in §3.7 it comes to constructing an $f$ for which $\mathcal{I}(f)$ has dimension 1 everywhere on $f(\partial\mathbb{H})$, that the same upper and lower bounds will hold even if we make (small) alterations in $X^{(N)}$. Let $M$ be an upper bound for $|g'(z)|$ on $\mathbb{H}$. From properties of $H_0^{-1}$ and $g$ it follows that there is a $\delta > 0$ such that $|g'(H_0^{-1}(w)) - b_0| < \frac{|b_0|}{4}$ for $|w| < \delta$. Let $T = \frac{24}{\pi\delta}$ and let $\eta$ be so small that
$$\tfrac{1}{\pi}(\tfrac{3}{T} + 2\eta) < \tfrac{\delta}{4}, \quad T\eta < \min\{1, \tfrac{|b_0|}{4M}\} \ \text{ and } \ \eta < \tfrac{|b_0|}{8M},$$
that is,

(27) $$\eta < \min\{\tfrac{\pi\delta}{24}, \tfrac{\pi\delta|b_0|}{96M}, \tfrac{|b_0|}{8M}\}$$

Note that $\delta$, $\eta$ and $T$ depend solely on $g$. If $Y \subset \mathbb{R}$ is a measurable set for which

(28) $$|\mathbf{P}(Y \setminus [0, b], z)| < \tfrac{\delta}{4}, \ \text{ for } z \in \mathbb{H} \cap \Delta(\tfrac{b}{2}, \tfrac{b}{2})$$

and the set $X = Y \cap [0, b]$ satisfies the hypothesis of the preceding lemma, then
$$|\mathbf{P}(Y, z)| < \tfrac{\delta}{2}, \ \text{ for } z \in F = [\tfrac{b}{N}, \tfrac{(N-1)b}{N}] + \tfrac{Tb\eta}{N}i,$$
so that with $G(z) = G(Y, z)$

$$|G((j+T\eta i)\tfrac{b}{N}) - G((k+T\eta i)\tfrac{b}{N})| \geq \tfrac{3}{4}|b_0||j-k|\tfrac{b}{N}.$$

(Note that we have used the fact that $F \subset \mathbb{H} \cap \Delta(\tfrac{b}{2}, \tfrac{b}{2})$, which follows by a trivial calculation since we have stipulated that $T\eta < 1$.) But then

$$\begin{aligned}|G(j\tfrac{b}{N}) - G(k\tfrac{b}{N})| &\geq \tfrac{3}{4}|b_0||j-k|\tfrac{b}{N} - 2M\tfrac{Tb\eta}{N} = (\tfrac{3}{4}|b_0||j-k| - 2MT\eta)\tfrac{b}{N} \\ &> (\tfrac{3}{4}|b_0||j-k| - \tfrac{|b_0|}{2})\tfrac{b}{N} \geq \tfrac{1}{4}|b_0||j-k|\tfrac{b}{N}.\end{aligned}$$

Since $\eta < \tfrac{|b_0|}{8M}$ we have (with $J = (-1, 1)$, as defined in §3.2) that

$$\begin{aligned}\operatorname{dist}(G(j\tfrac{b}{N} + \tfrac{\eta b}{2N}J), G(k\tfrac{b}{N} + \tfrac{\eta b}{2N}J)) &\geq \tfrac{1}{4}|b_0||j-k|\tfrac{b}{N} - 2M\tfrac{\eta b}{2N} > \tfrac{b_0}{8}|j-k|\tfrac{b}{N} \\ &> \tfrac{b_0}{8}\operatorname{dist}(j\tfrac{b}{N} + \tfrac{\eta b}{2N}J, k\tfrac{b}{N} + \tfrac{\eta b}{2N}J).\end{aligned}$$

If we let $b = \beta_1$ and $\epsilon \leq \tfrac{\eta}{2\pi}$ (so that $12\epsilon < \tfrac{\delta}{4}$ which allows us to apply (23) and therefore have (28) hold with $b = \beta_1$) we immediately conclude that $G(X^{(N)}, z)$ is bi-Lipschitz on $E^{(N)}$, that is, that for $N \geq 3$ and with the indicated values of $\epsilon$ and $b$,

(29) $\quad \tfrac{b_0}{8}|x_1 - x_2| \leq |G(X^{(N)}, x_1) - G(X^{(N)}, x_2)| \leq M|x_1 - x_2|, \ x_1, x_2 \in E^{(N)}.$

It therefore follows from the preservation of Hausdorff dimension under bi-Lipschitz mappings and the lower bound for $\dim(E^{(N)})$ obtained in §3.3 that

(30) $\quad \dim(G(X^{(N)}, E^{(N)})) = \dim(E^{(N)}) \geq d(N) = \tfrac{\log(N-1)}{\log N + \log(1+\log N) + \log(2/\epsilon\gamma_1)}.$

We will make use of the following observation in establishing that for the functions of the conclusion of Theorem A the set $\mathcal{I}(f)$ has dimension 1 locally in $f(\partial D)$.

**Observation 1**. It is clear form the foregoing considerations that for any fixed $N \geq 3$ every open interval $U$ of $\mathbb{R}$ has a subinterval $U'$ such that (29), with the constant $\tfrac{b_0}{8}$ replaced by any smaller number, and consequently (30), also hold for any $Y$ for which $X^{(N)}\Delta Y \subset U'$.

### 3.6. $\mathcal{I}(G^{(N)}) \supset G(X^{(N)}, E^{(N)})$.

Before beginning we note that the following weak version of Theorem A follows immediately from this inclusion together with (30): For any $\sigma < 1$ there is a mapping $f(z) = \int_0^z g_0'(\rho(\zeta))\,d\zeta$ for which $\dim(\mathcal{I}(f)) > \sigma$.

For $0 < c < 1$ let $D(c) = \overline{(\mathbb{H} + ci)} \cap \Delta$, for $0 \leq \beta < 1$ let $I(\beta) = J\setminus[-\beta, \beta]$, $H_\beta(z) = \mathbf{P}(I(\beta), z)$ and $g_\beta(z) = G(H_\beta, z)$, so that in particular $g_0 = g$. In addition, we denote by $\mathcal{U}_\beta$ the class of functions $u$ on $\mathbb{H}$ for which $u + H_\beta$ maps $\mathbb{H}$ into $\mathbb{S}$, and we denote by $\|u\|_A$ the sup-norm of $u$ on the set $A \subset \mathbb{H}$. By the normalization conditions (G*i*) and (G*ii*) of Lemma 6 there are positive constants $M > 1$ and $m$ for which $g'(\mathbb{H}) \subset R(m, M)$, so that $G'(H_\beta + u, z)$ satisfies the same inclusion for all $u \in \mathcal{U}_\beta$. Clearly, there exist $\beta_2, \xi \in (0, 1)$ which depend solely on $g$ such that for all $\beta \in [0, \beta_2]$ we have

$$g_\beta(z_\beta) = g_\beta(0) = 0 \text{ for some } z_\beta \in \mathbb{H} \text{ with } \Delta(z_\beta, \xi) \subset D(\xi).$$

It is a simply verified general fact that there exist positive constants $A$ and $B$ which depend only on $m$ and $M$ such that if $f'(\Delta) \subset R(m, M)$ then $f$ is univalent on $\Delta(A)$ and $f(\Delta(A)) \supset \Delta(f(0), B)$. Thus, there are positive constants $\rho < \frac{\xi}{2}$ and $r$, which depend only on $m$ and $M$ (and consequently only on $g$) such that for all $\beta \in [0, \beta_2]$ and all $u \in \mathcal{U}_\beta$, $G(H_\beta + u, z)$ is univalent on $\Delta(z_\beta, \rho)$ and

(31) $$\Delta(G(H_\beta + u, z_\beta), r) \subset G(H_\beta + u, \Delta(z_\beta, \rho)).$$

Note that $\Delta(z_\beta, \rho) \subset D(\xi)$. For any $c \in (0, 1)$ the set $H_0(D(c)) \subset \mathbb{S}$ is compact, so that simple considerations involving the continuous dependence with respect to $\beta$ of $H_\beta(K)$ for any fixed compact $K \subset \mathbb{H}$ show that there are positive numbers $\epsilon_1 = \epsilon_1(c)$, $M_1 = M_1(c)$, $\beta_3 = \beta_3(c) \leq \beta_2$ such that the derivative of $g'(H_0^{-1}(z))$ is bounded above by $M_1$ in the $\epsilon_1$-neighborhood of $\bigcup_{\beta \leq \beta_3} H_\beta(D(c)) \subset \mathbb{S}$ for all $\beta \leq \beta_3$ and $u \in U_\beta$. Thus, if $\beta \leq \beta_3$ and $u_\beta \in \mathcal{U}_\beta$ satisfies $\|u\|_{D(c)} \leq \epsilon_1$, then

$$|G'(H_\beta + u, z) - G'(H_\beta, z)| < M_1 \|u\|_{D(c)}, \quad z \in D(c).$$

Upon taking into account that $|G'(H_\beta + u, z)|$ is bounded above by $M$ on the segment $[0, ci]$ and that $[ci, z] \subset D(c)$ and $|z - ci| < 1$ for $z \in D(c)$, by integration over the contour $[0, ci] \cup [ci, z]$ we conclude that for such any $u$

$$|G(H_\beta + u, z) - G(H_\beta, z)| < 2Mc + M_1 \|u\|_{D(c)}, \quad z \in D(c).$$

Let $\eta$ be as in (27), $c_0 = \min\{\frac{r}{8M}, \frac{\xi}{4}\}$ and $\epsilon_0 = \min\{\frac{r}{4M_1}, \epsilon_1(c_0), \frac{\eta}{2\pi}\}$; $c_0$ and $\epsilon_0$ depend only on $g$. (We have included $\frac{\eta}{2\pi}$ in the minimum defining $\epsilon_0$ so that (30) holds.) Note that $c_0 < r$ since we have stipulated that $M > 1$ and that $\Delta(z_\beta, \rho) \subset D(c_0)$ since $\Delta(z_\beta, \rho) \subset D(\xi) \subset D(\frac{\xi}{4}) \subset D(c_0)$. Then for $\beta \leq \beta_3(c_0)$

$$|G(H_\beta + u, z) - G(H_\beta, z)| < \tfrac{r}{2}, z \in D(c_0),$$

for all $u \in \mathcal{U}_\beta$ for which $\|u\|_{D(c_0)} < \epsilon_0$. Since $G(H_\beta, z_\beta) = 0$, we have from this that for such $u$, $|G(H_\beta + u, z_\beta)| < \frac{r}{2}$. Also,

$$G(H_\beta + u, \overline{\mathbb{H} \cap \Delta(c_0)}) \subset \Delta(Mc_0) \subset \Delta(\tfrac{r}{2}),$$

so that by (31) we have that

$$G(H_\beta + u, \overline{\mathbb{H} \cap \Delta(c_0)}) \subset \Delta(\tfrac{r}{2}) \subset \Delta(G(H_\beta + u, z_\beta), r) \subset G(H_\beta + u, \Delta(z_\beta, \rho)).$$

Note that since $\Delta(z_\beta, \rho) \subset D(\xi)$ and $\overline{\Delta(c_0)} \subset \overline{\Delta(\xi/4)}$ and $\Delta(\xi/4) \cap D(\xi) = \emptyset$,

$$\overline{\mathbb{H} \cap \Delta(c_0)} \cap \Delta(z_\beta, \rho) = \emptyset.$$

We define

$$\beta_1 = \tfrac{1}{4}\min\{c_0, \beta_3(c_0)\}, \quad \epsilon = \tfrac{1}{4}\min\{\tfrac{\epsilon_0}{24}, \tfrac{1}{200}\}.$$

(The factors of $\tfrac{1}{4}$ are unnecessary at this point in that the relations (32) and (33), which will appear presently, hold without this factor; the reason for its inclusion will be explained below in Observation 2.) Here again $\beta_1$ and $\epsilon$ depend solely on $g$. As in §3.2, $\alpha = \frac{\epsilon}{N(1+\log N)}$ and for $\gamma_1$ we take any positive number less than $\epsilon\beta_1$ which is so small that if $U$ is any subset of $(-\tfrac{1}{2}\beta_1, \tfrac{1}{2}\beta_1)$ of measure at most $\gamma_1$, then

$$|\mathbf{P}(U, z)| < \epsilon \text{ on } \overline{\mathbb{H} \cap (\Delta \backslash \Delta(\beta_1))},$$

and finally, $\gamma = \tfrac{1}{2}\gamma_1\alpha$. We emphasize that $\epsilon$ *and* $\gamma_1$ *depend solely on* $g$ *and do not depend on* $N$. Also,

(32) $$\overline{\mathbb{H} \cap \Delta(2\beta_1)} \cap \Delta(z_{\beta_1}, \rho) = \emptyset$$

and

(33) $$G(H_\beta + u, \overline{\mathbb{H} \cap \Delta(2\beta_1)}) \subset G(H_\beta + u, \Delta(z_\beta, \rho)).$$

We can now use the bound (23) of §3.2 to show that

(34) $$\mathcal{I}(G^{(N)}) \supset G(X^{(N)}, E^{(N)}).$$

For the moment we abbreviate $X^{(N)}$ by $X$ and $G(X^{(N)}, z)$ by $G$ and examine the behavior of $G$ on the disk $\Delta(c(\mathbf{k}), \tfrac{\alpha}{2}\gamma^{l-1})$ where $l = L(\mathbf{k})$. We further abbreviate $c' = c(\mathbf{k})$ and $\gamma' = \tfrac{\alpha}{2}\gamma^{l-1}$. Let $z \in \Delta(c', \gamma')$. Now, for $\zeta \in c' + \overline{\mathbb{H} \cap \gamma'(\Delta\backslash\Delta(\beta_1))}$ we have by (23) of §3.2 that $\mathbf{P}(X, \zeta) = \mathbf{P}(\gamma' I + c', \zeta) + v(\zeta)$, where $\|v\|_{c' + \overline{\mathbb{H}\cap\gamma'(\Delta\backslash\Delta(\beta_1))}} < 12\epsilon \leq \tfrac{\epsilon_0}{2}$. But by (3) of §3.1, $\mathbf{P}(\gamma' I + c', \zeta) = \mathbf{P}(I, \tfrac{\zeta-c'}{\gamma'})$ so that if we substitute $\zeta = c' + \gamma'\omega$ we have

$$\begin{aligned}
G(z) - G(c') &= \int_{c'}^{z} g'(H_0^{-1}(\mathbf{P}(X, \zeta)))\, d\zeta \\
&= \gamma' \int_0^{(z-c')/\gamma'} g'(H_0^{-1}(\mathbf{P}(X, c' + \gamma'\omega)))\, d\omega \\
&= \gamma' \int_0^{(z-c')/\gamma'} g'(H_0^{-1}(\mathbf{P}(\gamma' I + c', c' + \gamma'\omega) + v(c' + \gamma'\omega)))\, d\omega
\end{aligned}$$

$$= \gamma' \int_0^{(z-c')/\gamma'} g'(H_0^{-1}(\mathbf{P}(I,\omega) + v(c' + \gamma'\omega))) \, d\omega.$$

Writing $v(c' + \gamma'\omega) = u(\omega)$, we have that $\|u\|_{\overline{\mathbb{H} \cap (\Delta \setminus \Delta(\beta_1))}} < \frac{\epsilon_0}{2}$, so that by (32) and (33) it follows that

(35) $$(\gamma'(\overline{\mathbb{H} \cap \Delta(2\beta_1)}) + c') \cap (\gamma'\Delta(z_{\beta_1}, \rho) + c') = \emptyset$$

(36) $$G(\gamma'\Delta(z_{\beta_1}, \rho) + c') \supset G(\gamma'(\overline{\mathbb{H} \cap \Delta(2\beta_1)}) + c').$$

From this and in light of the nesting properites of the intervals $J(\mathbf{k})$, $k \in \mathcal{K}$ together with (32) we conclude that (34) is true (and would remain true even without the factor of 2 appearing in these two equations).

**Observation 2**. This observation parallels Observation 1, which appears after (30). In light of the parenthetical comment made immediately after the definitions of $\beta_1$ and $\epsilon$ above it is clear that for any fixed $N \geq 3$, every open interval $U$ of $\mathbb{R}$ has a subinterval $U'$ such that (34) also holds if $X^{(N)}$ is replaced by any $Y$ for which $X^{(N)} \Delta Y \subset U'$. Indeed, since $E^{(N)}$ is nowhere dense, there is a point $p \in U \setminus E^{(N)}$, and it is clear $p \in \overline{J(\mathbf{k})}$ only for $\mathbf{k}$ in some finite (possibly empty) set $\mathcal{F}$ of nodes. In light of the continuous dependence of $G(X, z)$ on $X$ (that is, the fact that $G(X, z)$ will change by an arbitrarily small amount if $X$ is changed on a sufficiently small set) there is a neighborhood $V$ of $p$ such that $X^{(N)}$ can be changed arbitrarily on $V$ without altering the validity of (35) and (36) for any of the $\mathbf{k} \in \mathcal{F}$. Moreover, since there is a $\delta > 0$ such that $\text{dist}(p, \overline{J(\mathbf{k})}) \geq \delta$ for all $\mathbf{k} \notin \mathcal{F}$, it is clear that given any $\epsilon' > 0$ there is a $U' \subset V$ such that for any $Y \subset U'$, $|\mathbf{P}(Y,z)| < \epsilon'$ on all of the semi-disks $\Delta(c(\mathbf{l}), \frac{\alpha}{2}\gamma^{L(\mathbf{l})-1})$, $\mathbf{l} \notin \mathcal{F}$. From this together with the preceding sentence the desired conclusion follows.

**3.7 An $f$ for which $\mathcal{I}(f)$ has dimension 1 everywhere on $f(\mathbb{R})$.**

Let the sequence $\{p_n : n \geq 3\}$ be dense in $\mathbb{R}$ and and let $\sigma > 0$. The the construction involves the insertion into $\mathbb{R}$ of a copy $X_n = a_n X^{(n)} + c_n$, of $X^{(n)}$, successively for $n = 3, 4, \ldots$. Here, $0 < a_n$, $|c_n - p_n| \leq \frac{\delta}{2^n}$, from which the desired property will follow immediately. Because of the density of $\{p_n\}$, once we have inserted a copy $X_k$ of $X^{(k)}$, we will inevitably have to alter this copy by inserting copies of infinitely many $X^{(l)}$ with $l > k$ into the intervals making up $X_k$ and into the spaces between them. This must be done in such a way as not to alter the set of points of the corresponding copy $J_k$ of $(0, 1)$ whose images are assumed infinitely often (that is, the corresponding copy $E_k = a_k X^{(k)} + c_k$ of $E^{(k)}$). The process is nevertheless very straightforward and the existence of the desired mapping of $\mathbb{H}$ is an immediate consequence of the following construction together with the the facts, pointed out in Observations 1 and 2 at the ends of §3.5 and §3.6, respectively, to the effect that the bi-Lipschitz and infinite valence properties are unaffected by small changes in any of the $X^{(N)}$. To facilitate the discussion let $I' = J \setminus [-\frac{1}{2}, \frac{1}{2}]$, where as in §3.2, $J = (-1, 1)$, and for any $\rho > 0$ and $\tau \in \mathbb{R}$, let $I'(\rho, \tau) = \rho I' + \tau$. Let $Y$ be a subset of $\mathbb{R}$ whose complement contains the closed interval $\rho \overline{J} + \tau$. Since $\mathbf{P}(Y \cup I'(\rho,\tau), x)$ is pure imaginary on $\rho(-\frac{1}{2}, \frac{1}{2}) + \tau$ and tends to to $+\infty i$ and $-\infty i$ as $x$ tends to the left and right endpoints of this interval, there is a point $q \in \rho(-\frac{1}{2}, \frac{1}{2}) + \tau$, for which $\mathbf{P}(Y \cup I'(\rho,\tau), q) = 0$. Thus, for sufficiently small $\delta > 0$, $|\mathbf{P}(Y \cup I'(\rho,\tau), x)|$ will be arbitrarily small on $\overline{\mathbb{H} \cap (\Delta(\delta) + q)}$, so that for sufficiently small $\delta' > 0$, $\mathbf{P}(Y \cup I'(\rho,\tau) \cup (\delta' X^{(N)} + q), z)$ will have the same bi-Lipschitz and infinite valence properties on $\overline{\mathbb{H} \cap (\Delta(\delta') + \tau)}$ as $\mathbf{P}(X^{(N)}, z)$ does on $\overline{\mathbb{H} \cap \Delta(\frac{1}{2}, \frac{1}{2})}$.

Using this construction we proceed as follows. For $Y \subset \mathbb{R}$ let $\partial Y$ denote the boundary of $Y$ as a subset of $\mathbb{R}$, and let $Y_3 = X^{(3)}$. Assume that we have defined $Y_{N-1}$. Let $\tau_N$ be a point of $\mathbb{R} \setminus \partial Y_{N-1}$ for which $|\tau_N - p_N| \leq \frac{\sigma}{2^N}$. There are two cases.

Case 1. $\tau_N \notin Y_{N-1}$. In this case there is an open interval containing $\tau_N$ whose closure is disjoint from that of $Y_N$ and we define $Y_N = Y_{N-1} \cup (\delta_N X^{(N)} + \tau_N)$ with an appropriately small $\delta_N$.

Case 2. $\tau_N \in Y_{N-1}$. In this case there is an open interval $U$ containing $\tau_N$ whose closure is entirely contained in $Y_{N-1}$. We first remove a subinterval $U'$ of $U$ containing $\tau_N$ which is so small that its deletion will not disturb the bi-Lipschitz and infinite valence properties of $Y_{N-1}$. Then we proceed as Case 1.

It is easy to see that $G(Y_N, z)$ converges uniformly on bounded sets of $\overline{\mathbb{H}}$ to a function $G$ which has the desired properties.

**3.8 Smoothly bounded Jordan domains.**

Let $w$ map $\mathbb{H}$ one-to-one onto the given smoothly bounded Jordan domain $D$. Here $w'$ will be "very well-behaved" on $\mathbb{R} = \partial \mathbb{H}$. The function that will accomplish the desired mapping will be $f(z) = G_D(w^{-1}(z))$, where

$$G_D(z) = \int_0^z G'(\zeta) w'(\zeta) \, d\zeta,$$

$G(z)$ being the function constructed in the preceding section. Obviously,

$$f'(z) = G'_D(w^{-1}(z))(w^{-1})'(z) = G'(w^{-1}(z)) w'(w^{-1}(z))(w^{-1})'(z) = G'(w^{-1}(z)).$$

We also have that for each $Y_N$ constructed above the corresponding $G'(Y_N, z) = g'(H_0^{-1}(\mathbf{P}(Y_N, z)))$ and so is of the form $g'_0(\rho_1(z))$, where $\rho_1$ again maps $\mathbb{H}$ into itself and from this it follows that $f'$ has the desired form $g_0(\rho(z))$.

Finally, we have to verify that $\mathcal{I}(f)$ has the desired properties. This amounts to verifying that the $G_D$ has the desired properties on $\mathbb{H}$. This is simply a matter of correctly selecting the sequence $\{\delta_k\}$ in light of the following observation. Since $w'$ is continuous and nonvanishing on $\overline{\mathbb{H}}$ and $G'$ bounded there, for any $\epsilon'$ every point $x_0 \in \mathbb{R}$ has a neighborhood $U(x_0)$ in $\mathbb{C}$ such that for some $C \neq 0$,

$$|G_D(z_2) - G_D(z_1) - C(G(z_2) - G(z_1))| < \epsilon' |z_2 - z_1|, \text{ for } z_1, z_2 \in \overline{\mathbb{H} \cap U(x_0)}.$$

The desired conclusion is now obvious in light of and (35) and (36) (see the parenthetical comment in the sentence immediately following these tow relations) and the fact that in arbitrarily small neighborhoods of each point we have inserted copies of infinitely many of the $X^{(N)}$.

**4. Application to the Becker Criterion.** Let $q : \mathbb{H}/i \to \mathbb{C}$ be the extremal function defined in [BP], so that $f_1(z) = iq(z/i)$ displays the analogous extremal behavior for $\mathbb{H}$ (rather than for $\mathbb{H}/i$, which is the canonical half-plane of that paper). As shown there $|f_1''(z)/f_1'(z)| \leq 1/(2\Im\{z\})$ and for $\tau > 1$ arbitrarily close to 1 the function $f_\tau(z) = \int_0^z (f_1'(\zeta))^\tau \, d\zeta$ is not univalent. Obviously, $|f_\tau''(z)/f_\tau'(z)| \leq \tau/(2\Im\{z\})$ in $\mathbb{H}$. Consider $R_B = f_1'(\mathbb{H}) = q'(\mathbb{H}/i),$ which is a simply covered domain in $\mathbb{C}$. We shall make some comments about an important property of $R_B$ after presenting the following lemmas, in light of which Theorem C is an immediate corollary of Theorem A and Lemma 6.

**Lemma 8.** Let $p : \Delta \to \mathbb{H}$. If $f'(z) = f_\tau'(p(z))$, then $|f''(z)/f'(z)| \leq \frac{\tau}{1-|z|^2}$ in $\Delta$.

**Proof.** Let $0 < r < 1$. Let $u$ satisfy $u'(z) = f_1'(p(z))$. The set $p(r\Delta) \subset \mathbb{H}$ is compact, so that for some sufficiently large $K > 0$, $R_r = f_1'(p(r\Delta)) \subset f_1'(K\Delta + Ki)$, since the sets $K\Delta + Ki \subset \mathbb{H}$, $K > 0$ exhaust $\mathbb{H}$. Let $g(z) = u(rz)/r$, so that $g'(z) = u'(rz)$ and therefore $g'(\Delta) \subset R_r$. Now, $W_1(z) = f_1'(Kz + Ki)$ maps $\Delta$ one-to-one onto the domain $f_1'(K\Delta + Ki) \supset R_r$. Thus, there is a function $h : \Delta \to \Delta$ such that $g'(z) = W_1(h(z))$. We have

$$\begin{aligned}
|(\log(g'(z)))'| &= |W_1'(h(z))/W_1(z)||h'(z)| \\
&= K|f_1''(Kh(z) + Ki)/f_1'(Kh(z) + Ki)||h'(z)| \\
&\leq K|h'(z)|/(2\Im\{Kh(z) + Ki\}) = |h'(z)|/2\Im\{h(z) + i\} \\
&\leq |h'(z)|/(2(1 - |h(z)|)) = |h'(z)|\frac{1+|h(z)|}{2(1-|h(z)|^2)} \\
&\leq |h'(z)|/(1 - |h(z)|^2) \leq \frac{1-|h(z)|^2}{1-|z|^2}/(1 - |h(z)|^2) = \frac{1}{1-|z|^2},
\end{aligned}$$

the last inequality following from the Schwarz-Pick theorem. Since $r < 1$ is arbitrary we have that $|u''(z)/u'(z)| = |(\log(u'(z)))'| \leq \frac{1}{1-|z|^2}$. But

$$f'(z) = f_\tau'(p(z)) = (f_1'(p(z)))^\tau = (u'(z))^\tau,$$

so that $|f''(z)/f'(z)| = \tau|u''(z)/u'(z)| \leq \frac{\tau}{1-|z|^2}$. ∎

**Lemma 9.** Let $p : \mathbb{H} \to \mathbb{H}$. If $f'(z) = f_\tau'(p(z))$, then $|f''(z)/f'(z)| \leq \frac{\tau}{2\Im\{z\}}$ in $\mathbb{H}$.

**Proof.** This is an immediate consequence of the half-plane version of the Schwarz-Pick theorem which says that $|p'(z)| \leq \frac{\Im\{p(z)\}}{\Im\{z\}}$. ∎

We end with a few comments about first-order univalence criteria in general and about what we shall call the *first-order Becker criterion* in particular. A *first-order univalence criterion* for a domain $D$ is a condition of the form $f'(D) \subset R$ which implies that $f$ is univalent in $D$; such an $R$ is called a *domain of univalence* for $D$. In [G4] we defined such a criterion to be *weakly sharp* if for any $\epsilon > 0$ there is a nonunivalent function $f_\epsilon$ on $D$ for which $f'_\epsilon(D)$ is contained in the $\epsilon$-neighborhood of $R$, to be *strongly sharp* if for any domain $R'$ properly containing $R$ there are is a nonunivalent function $f$ for which $f'(D) \subset R'$, and, finally, to be *very strongly sharp* if for any such $R'$ there are functions $f$ for which $\mathcal{I}(f)$ is not empty. In light of Theorem A it is clear that it would be appropriate to alter the last definition by changing the words "for which $\mathcal{I}(f)$ is not empty" to "for which $\mathcal{I}(f)$ has dimension at least 1" or to introduce some other term, such as "super-strongly sharp" to describe the case in which this last condition holds. As in [G4, p.174], let $\mathfrak{R}$ denote the class of domains $R$ of the form $R = S_1 \setminus \overline{S_2}$, where $S_1$ and $S_2$ are analytically bounded strictly star-shaped (with respect to 0) domains with $\overline{S_2} \subset S_1$. Then in light of the corollary to Theorem A, Theorem 6 of that paper can be changed to read

**Theorem.** *If $f'(\mathbb{H}) \subset R \in \mathfrak{R}$ is a weakly sharp univalence criterion and $R$ is a domain of univalence for $D$, then for any domain $R'$ properly containing $R$ there is a function $f$ for which $f'(D) \subset R'$ and for which $\dim(\Delta(z_0, r) \cap \mathcal{I}(f)) = 1$ for all $r > 0$ and all $z_0 \in f(\partial D)$.*

We will not go into a discussion of the exact family of simply connected domains $D$ to which this applies; it certainly contains all smoothly bounded Jordan domains, though.

If one defines a first-order criterion $f'(D) \subset R$ to be *very weakly sharp* if there is a sequence of domains $\{R_n\}$ none of which is a domain of univalence for $D$ and such that $R = \cap R_n$, then it follows from [BP] that the *first-order Becker criterion*, that is, the condition $f'(D) \subset R_B$, is a very weakly sharp criterion both for $D = \Delta$ (in light of Lemma 8) and for $D = \mathbb{H}$. Because of the square root singularities of the extremal function $f_1$ at each of the two points $\pm \pi \in \partial \mathbb{H}$ for which $f_1(\pi) = f_1(-\pi)$, the treatment of [G4], which is based on an analysis of the signs of the appropriate variation kernel does not go through for this criterion. We believe, however, that this criterion is actually weakly sharp but *not* strongly sharp. We will not go into details of the motivation for this conjecture beyond saying that it is based on a calculation of the variation kernel corresponding to $R_B$, which, in turn, is based on simple observations about the shape of this domain. To my knowledge, the criteria $f'(\Delta) \subset R_B$, $f'(\mathbb{H}) \subset R_B$ are the *only* sharp first-order criteria which are neither trivial (like the Noshiro-Warschawski criterion) nor ones for which the $R$ is some very thin domain (such as the those derived in [G2] and [G3]). In spite of this, I do not believe that this first-order version of the Becker criterion has been given explicit mention anywhere in the literature.

**Acknowledgement.** The author would like to thank Víctor Castillo of the Catholic University of Chile for carefully reading and correcting a preliminary version of this paper.

REFERENCES

[A] Aksent'ev, L. A.*Multivalent functions from extended Becker and Nehari classes and their hydromechanical interpretation*. Russian Math. (Iz. VUZ) 43 (1999), no. 6,

1–12.

[AK] Avhadiev, F. G. and Kayumov, I. R. *Admissible functionals and infinite-valent functions*, Complex Variables Theory Appl. 38 (1999), no. 1, 35–45.

[BP] Becker, J. and Pommerenke, Ch., *Schlichtheitskriterien und Jordangebiete*, J. Reine Angew. Math. 354 (1984), 74–94.

[F] Falconer, K., *Fractal Geometry, Mathematical Foundations and Applications*, John Wiley & Sons, Chichester, 1990

[G1] Gevirtz, Julian, *An upper bound for the John constant*, Proc. Amer. Math. Soc. 83 (1981), no. 3, 476–478.

[G2] Gevirtz, Julian, *On sharp first-order univalence criteria*, Indiana Univ. Math. J. 36 (1987), no. 1, 189–209.

[G3] Gevirtz, Julian, *On the difference quotients of an analytic function*, Ann. Acad. Sci. Fenn. Ser. A I Math. 12 (1987), no. 2, 237–259

[G4] Gevirtz, Julian, *The theory of sharp first-order univalence criteria*, J. Anal. Math. 64 (1994), 173–202.

[P] Pommerenke, Ch., *Boundary Behaviour of Conformal Maps*, Grundlehren der Mathematischen Wissenschaften 299, Springer-Verlag, Berlin, 1992.